\documentclass[letterpaper, paper,11pt]{AAS}		

\usepackage{bm}
\usepackage{amsmath}
\usepackage[colorlinks=true, pdfstartview=FitV, linkcolor=black, citecolor= black, urlcolor= black]{hyperref}
\usepackage{overcite}
\usepackage{footnpag}
\usepackage{placeins}
\usepackage[colorlinks=true, pdfstartview=FitV, linkcolor=black, citecolor= black, urlcolor= black]{hyperref}
\usepackage{breqn}
\usepackage{float}
\usepackage{graphicx}
\graphicspath{{images/}{../images/}}
\usepackage{subfiles}
\usepackage{subcaption}
\usepackage{comment}
\usepackage{algorithm}
\usepackage{algpseudocode}
\usepackage{graphics} 
\usepackage{epsfig} 
\usepackage{mathptmx} 
\usepackage{times} 
\usepackage{amsmath} 
\usepackage{amssymb}  
\usepackage{xcolor,colortbl}
\usepackage{subcaption}
\usepackage{comment}

\newcommand{\x}{\textbf{x} }

\newcommand{\Q}{{\textbf{Q}} }

\newcommand{\be}{\begin{equation}}
\newcommand{\ee}{\end{equation}}

\usepackage{enumitem}
\usepackage{multirow,array}
\usepackage{algorithm,algcompatible,algorithmicx}
\algnewcommand\INPUT{\item[\textbf{Input:}]}%
\algnewcommand\OUTPUT{\item[\textbf{Output:}]}%

\usepackage{cite}
\usepackage{lipsum,booktabs}
\renewcommand{\u}{\mathbf{u}}

\newcommand{\reffig}[1]{Fig. \ref{#1}}

\newcolumntype{L}{>{$}l<{$}} 
\newcolumntype{C}{>{$}c<{$}}

\PaperNumber{25-291}
\begin{document}

\title{A Hamilton-Jacobi Approach for Nonlinear Model Predictive Control in Applications with Navigational Uncertainty}

\author{Amit Jain\thanks{Postdoctoral Associate, MIT Department of Aeronautics and Astronautics},  
\ Roshan T. Eapen\thanks{Assistant Professor, Department of Aerospace Engineering, Pennsylvania State University, University Park, PA-16802}
\ and Puneet Singla\thanks{Harry and Arlene Schell Professor of Engineering, Department of Aerospace Engineering, Pennsylvania State University, University Park, PA-16802.}
}

\maketitle{}

\begin{abstract}
This paper introduces a novel methodology that leverages the Hamilton-Jacobi solution to enhance non-linear model predictive control (MPC) in scenarios affected by navigational uncertainty.  Using Hamilton-Jacobi-Theoretic approach, a methodology to improve trajectory tracking accuracy among uncertainties and non-linearities is formulated. This paper seeks to overcome the challenge of real-time computation of optimal control solutions for Model Predictive Control applications by leveraging the Hamilton-Jacobi solution in the vicinity of a nominal trajectory. The efficacy of the proposed methodology is validated within a chaotic system of the planar circular restricted three-body problem.
\end{abstract}

The precise tracking of trajectories in nonlinear systems is a fundamental focus within control theory. In many real-world applications, such as robotics, autonomous vehicles, and aerospace, nominal trajectories are designed and optimized to satisfy only scientific requirements as well as to comply with system constraints. However, the nominal path will unlikely be followed by the spacecraft in real-life scenarios due to uncertainty in the dynamic model (e.g. gravitational parameters, environmental effects), navigational errors (e.g. imperfect state knowledge or approximations in the measurement models), command actuation (i.e. thrust magnitude, pointing angle errors), etc. Correction maneuvers are then needed to compensate for these deviations which typically involve solving another optimal control problem. Furthermore, the high nonlinearity of the vehicle dynamics makes it difficult for the trajectory tracking controllers to ensure tracking accuracy at all times. Therefore, a good trajectory tracking control system should fully consider the nominal trajectory design and the nonlinearities in the dynamics.

In tackling the trajectory tracking challenge, extensive research has focused on the selection of appropriate control algorithms, including the linear quadratic regulator (LQR) \cite{yue2016zero,park2021experimental}, fuzzy control \cite{yang2017trajectory,xiong2010intelligent}, sliding mode control (SMC) \cite{hwang2017path,bei2022research}, and proportional-integral-derivative (PID) control \cite{marino2009nested}. Model predictive control (MPC) methods have gained widespread adoption in recent years for trajectory tracking control. Characterized by continual optimization and consideration of input and output constraints, MPC maximizes performance criteria while ensuring system constraints remain unviolated \cite{falcone2008mpc}. Significant advancements have been made in MPC controller formulation, including constraint incorporation \cite{brown2017safe}, adaptive adjustment of weight matrices in the target function to enhance trajectory tracking accuracy \cite{zhang2019electrical, tan2018mpc}, utilization of dead-band penalty functions to address inequality constraints and ensure smoothness \cite{guo2020computationally}, among others. To enhance the controller's predictive capabilities, real-time predictive vehicle control has been explored utilizing Pontryagin's optimization method \cite{banginwar2022autonomous}.

Despite its great potential, the pervasive use of MPC in industrial applications is somewhat limited by the ability to solve optimal control problems in real-time. This work seeks to alleviate the challenges of real-time computation of optimal control solutions for MPC application through the use of the Hamilton-Jacobi solution in the vicinity of a nominal trajectory. The Hamilton-Jacobi solution provides a judicious coordinate transformation that encodes the information of all optimal control trajectories between two points in a predefined domain \cite{jain2023sparse, eapen2021Diss}. In this work, the domain is defined in the neighborhood of a nominal trajectory.  The optimal control utilizes the current state information as predicted by the dynamical system model and is modestly robust to model errors and state uncertainties.

The proposed approach utilizes the HJ equation to solve for generating functions for the underlying Hamiltonian system in state and co-states. The generating functions provide a map between the current state and co-state to its value at the initial time. Numerous strategies have been proposed for solving the HJ equation within the context of optimal control \cite{achdou2013hamilton,guibout2004solving, park2004solutions,vaidya2022spectral}.
 The special nature of the relationship between the value function and the generating function is clarified in \cite{guibout2004solving, park2004solutions}, where a family of value functions can be derived from a single generating function. Thus, solving the HJ equation for a generating function achieves a family of optimal feedback control profiles as an explicit function of the boundary conditions. Eapen et al. \cite{eapen2021semi,roshanEapen2021} have investigated this property of  Hamiltonian dynamical systems in the context of the optimal feedback control problem. By connecting the value function to the Hamilton-Jacobi generating function, a systematic way to evaluate the optimal feedback control and cost function while still satisfying the general boundary conditions was obtained. Recently Jain et al. \cite{jain2023sparse} have utilized the Conjugate Unscented Transformation (CUT) based sparse-collocation method to solve the optimal feedback control problem.

The remainder of the paper is structured as follows: initially, the formulation of non-linear trajectory tracking is presented, followed by a concise introduction to the Hamilton-Jacobi equation and generating functions. Subsequently, the paper details the proposed solution methodology, employing the CUT-based sparse-collocation method to determine the trajectory tracking solution. Lastly, a numerical solution is showcased for the planar circular restricted three-body problem.



\section{Problem Formulation and Hamilton-Jacobi Preliminaries}

The primary objective of this research is to develop a numerical framework for solving the trajectory tracking problem in nonlinear dynamical systems. This framework is based on the Hamilton-Jacobi (HJ) theory, which provides a powerful tool for deriving the optimal control laws and trajectories for such systems. The problem formulation begins with the identification of an optimal trajectory that transitions from an initial state $\x_0$ to a final state  $\x_f$ within a specified time $t_f$. In this work, the Planar Circular  Restricted Three-Body Problem (PCR3BP) is used to demonstrate the efficacy of the MPC formulation with the HJ theory. The PCR3BP describes the motion of a particle subjected to the gravitational attraction of two bodies. The PCR3BP equilibrium points ($L_1$ - $L_5$) host a variety of periodic orbits in their vicinity, and the motions near them are known to be chaotic. This work considers a transfer between two planar periodic orbits at $L_1$, and $L_2$, also called Lyapunov orbits. The following vector differential equation governs the PCR3BP dynamics \cite{eapen2022momentum}:
\begin{equation}
    \ddot{\mathbf{q}}+2 \omega \times \dot{\mathbf{q}}=\frac{\partial \Omega}{\partial \mathbf{q}} + \mathbf{u}
\end{equation}
where $\boldsymbol{\omega}=[0,0,1]^T$, $\mathbf{q}=[x,y]^T$, and $\mathbf{u} = [u_x, u_y]$. The pseudo-potential $\Omega$ for PCR3BP is expressed as
$ \frac{1}{2}\left(x^2+y^2\right)+\frac{1-\mu}{r_1}+\frac{\mu}{r_2}$, where $r_1$ and $r_2$ denote the distances from the third body to the first and second primaries, respectively. These distances are defined as
$r_1=\left((x+\mu)^2+y^2\right)^{1 / 2}$ and $r_2=\left((x-1+\mu)^2+y^2\right)^{1 / 2}$.  

\begin{figure}[htb]
    \includegraphics[width=\linewidth]{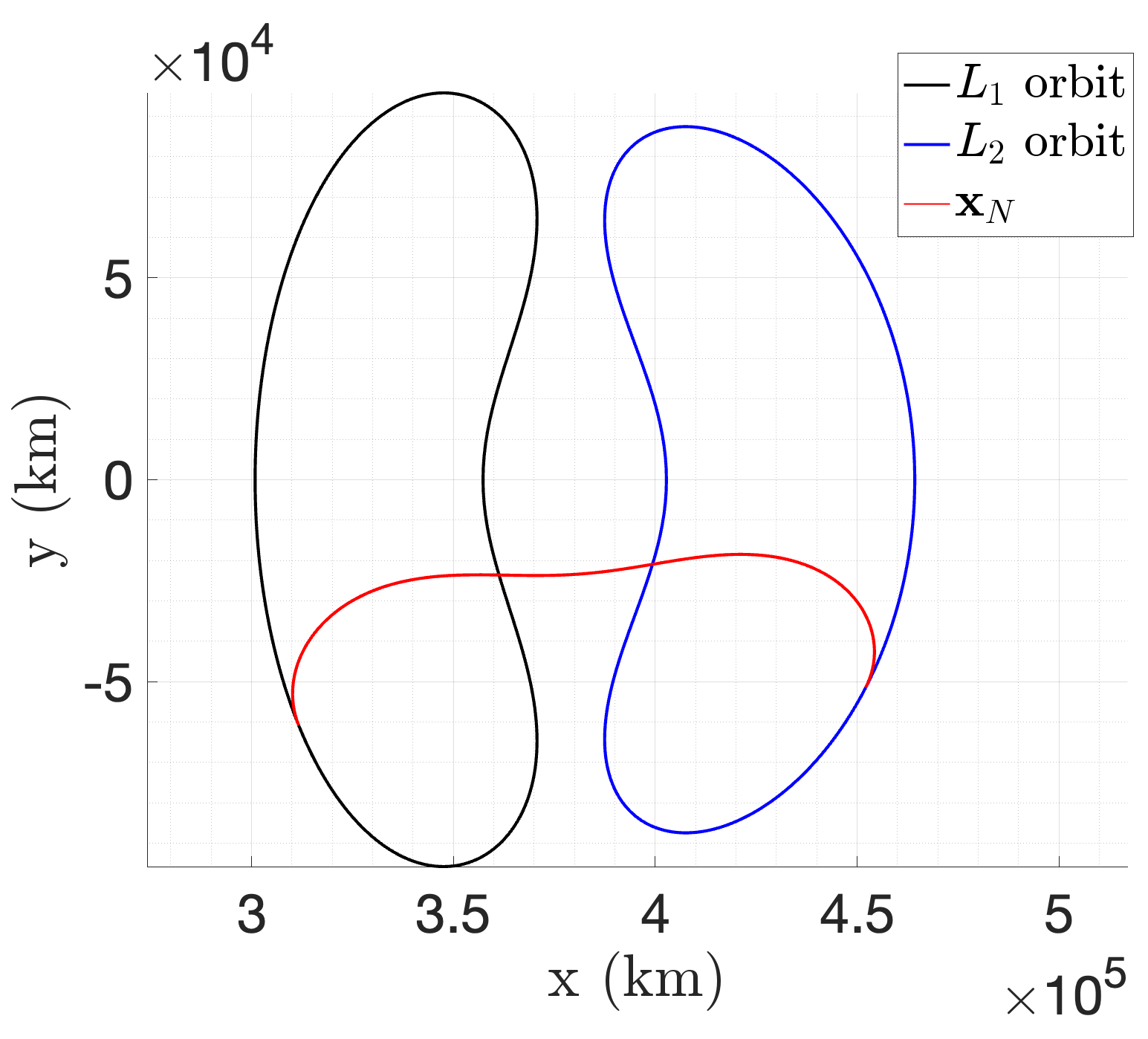}
    \caption{Nominal trajectory from $L_1$ to $L_2$ Lyapunov orbit}
    \label{fig:CR3BP_nominal}
\end{figure}

\reffig{fig:CR3BP_nominal} shows the nominal transfer trajectory from $\x_0$ in $L_1$  to $\x_F$ in $L_2$ Lyapunov orbit achieved by minimizing the cost function as shown in Table \ref{tab:CR3BP_OCP}. The generalized optimal control problem can be formulated by augmenting the boundary condition on the final state using the Lagrange multiplier, leading to the following formulation:
\begin{equation}\label{eq:OCP_track}
\min _{\mathbf{u}(t)} \quad   J = \phi({\mathbf{x}(t_f)}) +  \boldsymbol{\nu}^T \boldsymbol  \psi({\mathbf{x}(t_f)}) + \int_0^{t_f}\left(\mathbf{x}^T \Q \mathbf{x}+\mathbf{u}^T \mathbf{R u}\right) dt \\ 
\end{equation}
By solving the above optimal control problem, the nominal state trajectory $\mathbf{x}_N(t)$ and the corresponding control trajectory $\mathbf{u}_N(t)$ can be obtained over the time interval $t \in\left[t_0, t_f\right]$. 

\begin{figure*} 
\centering
\begin{subfigure}{\textwidth}
    \includegraphics[width=0.95\textwidth]{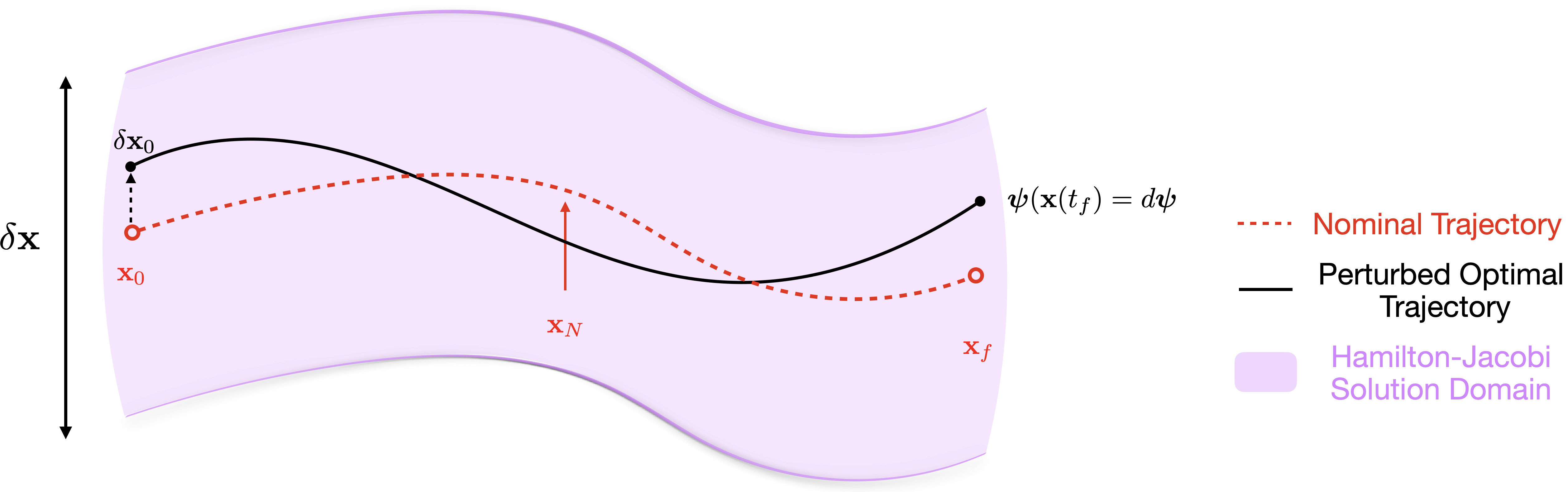}  
    \caption{Trajectory Tracking due to Initial Perturbation}
    \label{fig:trajtracking_a}
\end{subfigure}
  \hfill 
  \vspace{0.1in}
\begin{subfigure}{\textwidth}
  \centering
    \includegraphics[width=0.95\textwidth]{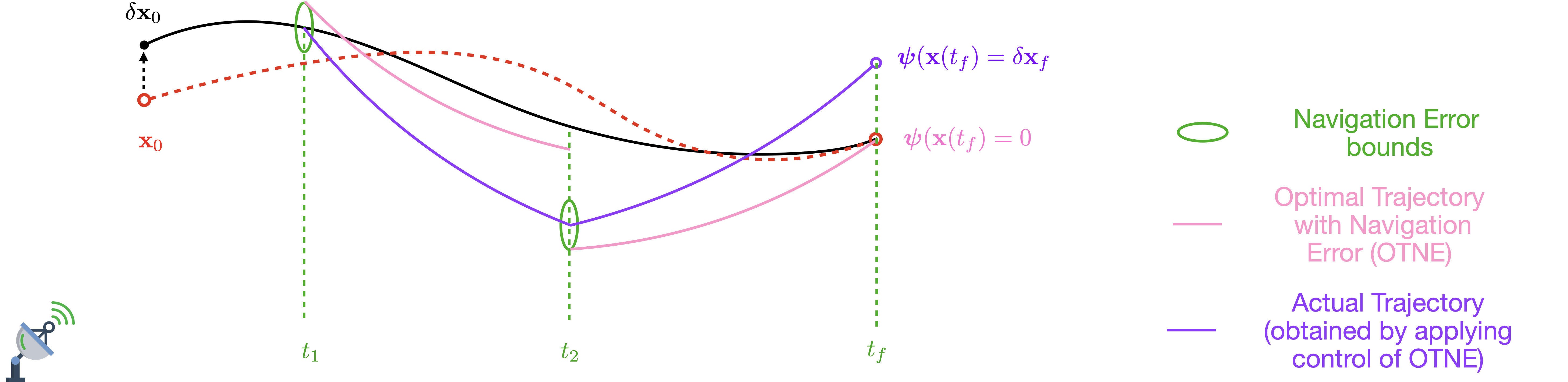} 
    \caption{Trajectory Tracking due to Initial Perturbation and Navigation Errors}
    \label{fig:trajtracking_b}
\end{subfigure}
\caption{Trajectory Tracking due to Perturbations and Errors}\label{fig:trajtracking}
\end{figure*}

Now that the nominal trajectory is designed, this paper will consider the case of an initially perturbed trajectory with navigational uncertainties. To accurately follow the nonlinear trajectory, establishing a nonlinear guidance problem in the proximity of the nominal trajectory is essential. Examining slight deviations from the nominal path caused by small perturbations in the initial state $\delta \mathbf{x}_0$ and terminal conditions $\delta \boldsymbol \psi$, as depicted in \reffig{fig:trajtracking_a}. The new cost function to be minimized can be expressed as:

\begin{subequations}\label{eq:J_secondorder}
\begin{align}
\min _{\mathbf{u}(t)} \quad  \Bar{J}  = & \frac{1}{2}\left[\delta \x^T\left(\phi_{\x \x}+\left(\boldsymbol \nu^T \boldsymbol \psi_\x\right)_\x\right) \delta \x\right]_{t=t_f}   +\frac{1}{2} \int_{0}^{t_f} \left(\delta \mathbf{x}^T \Q \delta \mathbf{x}+\delta \mathbf{u}^T \mathbf{R} \delta \mathbf{u}\right) dt  \\ \vspace{0.5in}
& \qquad  \text{subject to: }  \qquad  \delta \dot{\mathbf{x}}=\mathbf{f}(\mathbf{x}, \mathbf{u}) - \dot{\mathbf{x}}_N   
\end{align}
\end{subequations}
where $\delta \x = \x - \x_N$ and $\delta \u = \u - \u_N$. The Hamilton-Jacobi formulation is employed to solve the above non-linear trajectory tracking problem around the nominal trajectory, as illustrated in \reffig{fig:trajtracking_a}. The advantage of utilizing the HJ solution for computing the solution lies in its capacity to encapsulate the information of adjacent trajectories within the computed domain of interest $\delta \x$. This feature is further exemplified in \reffig{fig:trajtracking_b}, where the navigational uncertainty leads to an error in detecting the vehicle's actual position. The optimal trajectory with navigational error (OTNE) will be able to compute the control law that steers the state at each waypoint(shown on the boundary of the green error ellipse) to the desired final state where $\boldsymbol{\psi}(\mathbf{x}(t_f))=0$. However, if the actual position of the vehicle is at the center, the computed control law will take the vehicle to the undesired state (shown in purple color), where $\boldsymbol{\psi}(\mathbf{x}(t_f))=\delta \mathbf{x}_f$.

\begin{table}[htb!]
\centering
\caption{Types of Generating functions}
\setlength{\tabcolsep}{10pt}
\begin{tabular}{|C|L|C|}
\hline
\textbf{Type} & \multicolumn{1}{c|}{\textbf{Canonical Relationship}}  &\multicolumn{1}{c|}{  \textbf{Hamilton-Jacobi (HJ) equation} } \\ \hline  
& &   \\[0.01mm] 
F_1(\delta  \mathbf{x} , \delta  \mathbf{x}_0, t, t_0) & 
\displaystyle  \delta  \boldsymbol{\lambda}  =\frac{\partial F_1}{\partial \delta  \mathbf{x} },  \delta  \boldsymbol{\lambda}_0=-\frac{\partial F_1}{\partial \delta  \mathbf{x} _0} & 
\displaystyle  \frac{\partial F_1}{\partial t}+H\left(\delta  \mathbf{x} , \frac{\partial F_1}{\partial \delta  \mathbf{x} }, t\right)=0  \\[3mm]
\hline  
& &   \\[0.01mm] 
F_2(\delta  \mathbf{x} , \boldsymbol{\lambda_0}, t, t_0) & 
\displaystyle  \delta  \boldsymbol{\lambda}  =\frac{\partial F_2}{\partial \delta  \mathbf{x} } 
, 
\displaystyle \delta  \mathbf{x} _0=\frac{\partial F_2}{\partial \delta  \boldsymbol{\lambda}_0} &  
\displaystyle \frac{\partial F_2}{\partial t}+H\left(\delta  \mathbf{x} , \frac{\partial F_2}{\partial \delta  \mathbf{x} }, t\right)=0    \\[3mm]
\hline  
& &   \\[0.01mm] 
F_3(\delta  \boldsymbol{\lambda}, \mathbf{x_0}, t, t_0) & 
\displaystyle \delta  \mathbf{x}    =-\frac{\partial F_3}{\partial \delta  \boldsymbol{\lambda}},   \delta  \boldsymbol{\lambda}_0=-\frac{\partial F_3}{\partial \delta  \mathbf{x} _0} 
&
\displaystyle  \frac{\partial F_3}{\partial t}+H\left(-\frac{\partial F_3}{\partial \delta  \boldsymbol{\lambda}}, \delta  \boldsymbol{\lambda}, t\right)=0   \\[4mm] 
\hline 
& &   \\[0.01mm] 
F_4(\delta  \boldsymbol{\lambda}, \boldsymbol{\lambda_0}, t, t_0) & 
\displaystyle \delta  \mathbf{x}  =\frac{\partial F_4}{\partial \delta  \boldsymbol{\lambda}} ,
\delta  \mathbf{x} _0=-\frac{\partial F_4}{\partial \delta  \boldsymbol{\lambda}_0}
&
\displaystyle \frac{\partial F_4}{\partial t}+H\left(-\frac{\partial F_4}{\partial \delta  \boldsymbol{\lambda}}, \delta  \boldsymbol{\lambda}, t\right)=0 \\[5mm] \hline
\end{tabular}
\label{tab:GenFunc}
\end{table}

There exist four possible forms of the generating function which can generate the given canonical transformation to the new variable space ($\delta \mathbf{x}_0, \delta \boldsymbol{\lambda_0})$ and satisfy the HJ equation as shown in Table \ref{tab:GenFunc}\cite{jain2023sparse}. Each of these generating functions provides a different perspective on the system's dynamics and offers a unique approach to solving the optimal trajectory problem. The generating function of type-2, $F_2$, is particularly advantageous for problems where the initial and terminal perturbed states are explicitly given, as it directly relates the initial perturbed state and final perturbed costate, thereby solving the boundary-value problem.

The type-2 generating function $F_2$ is thus utilized to solve the trajectory tracking problem as summarized below:
\begin{enumerate}
    \item The Hamiltonian in $(\delta\mathbf{x}, \delta\boldsymbol\lambda, \delta\mathbf{u}, t $) is transferred to a new Hamiltonian through a canonical coordinate transformation:
    \begin{equation}\label{eq:Transformation_HJ} 
        \left(\delta\mathbf{x}(t), \delta\boldsymbol{\lambda}(t)\right) \rightarrow\left(\delta\mathbf{x}_{0}, \delta\boldsymbol{\lambda}_{0}\right)
    \end{equation}
    
    \item Using  $F = F_2(\delta\mathbf{x}, \delta\boldsymbol\lambda_0, t)$, the HJ equation and the corresponding relations:
       \begin{subequations}\label{eq:F2_lambda_x0}  
        \begin{align}
        \frac{\partial F_2\left(\delta\mathbf{x}, \delta\boldsymbol{\lambda}_0, t \right)}{\partial t}+ & H\left(\delta\mathbf{x}, \frac{\partial F_2\left(\delta\mathbf{x}, \delta\boldsymbol{\lambda}_0, t\right)}{\partial \delta\mathbf{x}}, t\right)=0 \label{eq:F2_HJE} \\
        & \delta\boldsymbol{\lambda}  =\frac{\partial F_2\left(\delta\mathbf{x}, \delta\boldsymbol{\lambda}_0, t \right)}{\partial \delta\mathbf{x}}  \label{eq:F2_lambda} \\
        & \delta\mathbf{x}_0=\frac{\partial F_2\left(\delta \mathbf{x}, \delta \boldsymbol{\lambda}_0, t \right)}{\partial \delta \boldsymbol{\lambda}_0} \label{eq:F2_x0} 
        \end{align}
        \end{subequations}
    
    \item Approximate $F$ as a polynomial: 
    \begin{subequations}\label{eq:F2_lambda_0_lambda_f_a}  
        \begin{align}
        F = F_2(\delta \mathbf{x}, \delta \boldsymbol\lambda_0, t_f) = F_2(\delta \mathbf{x}_f, \delta \boldsymbol\lambda_0) \\
        \delta \mathbf{x}_0 = \frac{\partial \delta F_2(\mathbf{x}_f, \delta \boldsymbol\lambda_0)}{\partial {\delta \boldsymbol{\lambda}_0}} = \kappa(\delta \mathbf{x}_f,\delta \boldsymbol{\lambda}_0)  \\
        \text{Rewriting: \qquad }  \delta \boldsymbol{\lambda}_0 = \bar{\kappa}(\delta \mathbf{x}_0, \delta \mathbf{x}_f)  \qquad \qquad \label{eq:lambda0(x0,xf)}
        \end{align}
        \end{subequations}

    \item So, solving for $F$ provides $\delta \boldsymbol{\lambda}_0$ as a function of ($\delta \mathbf{x}_0, \delta \mathbf{x}_f)$ i.e. any trajectory can be solved by finding $F$.
\end{enumerate}

The following section will outline the methodology for solving the HJ equation and obtain an expression for $F$ for optimal trajectory tracking. 

\section{Solution Methodology for MPC Trajectory Tracking using HJ Theory}

This section presents the methodology for trajectory tracking utilizing the type-2 generating function. The section begins by outlining the development of collocation equations, the selection of collocation points, and the optimal selection of basis functions.  These steps are crucial in formulating an effective approach for trajectory tracking, allowing us to achieve accurate and efficient results.

\subsection{Development of Collocation Equations} \label{CollocationEquations}
This problem aims to solve the HJ equation given in \eqref{eq:F2_HJE}. Following the definition of the approximate generating function, the solution to the HJ equation is assumed to be of the form:
\begin{equation}\label{eq:S(Q,Y)_track}
F_2(\delta  \mathbf{x}, \delta  \boldsymbol{\lambda}_0, t)=\sum_{j=1}^{m} \phi_{j}(\delta \mathbf{x}, \delta  \boldsymbol{\lambda}_0)  c_{j}(t) 
 =  \Phi^T(\delta  \mathbf{x}, \delta  \boldsymbol{\lambda}_0) \mathbf{c}(t)
\end{equation}
where  $\phi(\delta \mathbf{x}, \delta \boldsymbol{\lambda}_0) \in \mathbb{R}^m$ is a vector of the basis function and is assumed to have at least continuous first-order derivatives, while $\mathbf{c}(t) \in \mathbb{R}^m$ is a vector of time-varying coefficients. Substituting \eqref{eq:S(Q,Y)_track} in \eqref{eq:F2_HJE} and expanding the time derivatives as finite differences, the HJ equation can be written as:
\begin{subequations}\label{eq:FiniteDiff}
\begin{align}
    \Phi^T(\delta \mathbf{x}, \delta \boldsymbol{\lambda}_0) \dot{\mathbf{c}}_{i}(t)  & = - H\left(\delta \mathbf{x}, \frac{\partial F_{2}}{\partial (\delta \mathbf{x})}, t\right) \\
    \Phi^T(\delta\mathbf{x}, \delta\boldsymbol{\lambda}_0)  \frac{\mathbf{c}_{k+1} -\mathbf{c}_{k}}{dt}  &= -H\left(\delta\mathbf{x}, \frac{\partial F_{2}}{\partial (\delta\mathbf{x})}, t_k \right) 
 \\
    \Phi^T(\delta\mathbf{x}, \delta\boldsymbol{\lambda}_0)  \delta \mathbf{c}_{k+1}  &= -H_k dt \label{eq:delta_ck+1 Phi = K}  \\
    \mathbf{A} \delta \mathbf{c}_{k+1} & = \mathbf{b} \label{eq:FiniteDiff_d}
\end{align}
\end{subequations}
where $\mathbf{c}_{k}$ are the coefficients at time $t_k$ and $\delta \mathbf{c}_{k+1}=\mathbf{c}_{k+1} - \mathbf{c}_{k}$ are the departure coefficients at time $t_{k+1}$. The residual error from \eqref{eq:FiniteDiff_d} is given as:
\begin{subequations}
\begin{align}
\mathbf{e}(\boldsymbol{\zeta})  &= \mathbf{A} \delta \mathbf{c}_{k+1} - \mathbf{b} \label{eq:error_A_B}\\
\mathbf{A}^T_i  =\Phi^i(\boldsymbol{\zeta}), \quad
& \mathbf{b}_i  = -H\left(\delta\mathbf{x}, \frac{\partial F_{2}}{\partial (\delta\mathbf{x})}, t_k \right)  dt \label{eq:A_B}
\end{align}
\end{subequations}
where $\boldsymbol{\zeta} = [\delta\mathbf{x}, \delta\boldsymbol{\lambda}_0]$.
In the collocation method, the error is projected onto a series of delta functions centered at chosen collocation points, resulting in a residual error being zero at the collocation points. Assuming there are total $N$ collocation points, leading to a system of $N$ equations in $m$ unknowns to exactly solve the HJ equation at prescribed points, $\boldsymbol{\zeta}_i$: \vspace{-0.05in}
\begin{equation} \vspace{-0.07in}
\int e(\boldsymbol{\zeta}) \delta\left(\boldsymbol{\zeta}-\boldsymbol{\zeta}_i\right) d \boldsymbol{\zeta}=0 \rightarrow e\left(\boldsymbol{\zeta}_i\right)=0, \quad i=1,2, \ldots N
\end{equation}
where $\boldsymbol{\zeta}_i$ are the chosen collocation points. The selection of the collocation points is crucial in obtaining a well-conditioned system of equations for the unknown coefficients. Therefore, an efficient numerical sampling method known as the Conjugate Unscented Transformation (CUT) method is utilized to generate collocation points in $8$-dimensional space, providing a lower number of points compared to traditional quadrature methods like Gauss quadrature and sparse grid~\cite{adurthi2018conjugate}.


\subsection{Optimal Selection of Basis Functions}
The solution of Eq. \eqref{eq:error_A_B} can be calculated by minimizing the weighted two-norm error, which seeks to find the best-fit solution for the given collocation points:  
\begin{equation} 
\begin{aligned}
\delta \mathbf{c}_{{k+1}_{l_2}} = & \min_{\delta \mathbf{c}_{{k+1}}}  \left\|\mathbf{W}( \mathbf{A} \delta \mathbf{c}_{k+1} - \mathbf{b} ) \right\|_{2} \\
\end{aligned}
\end{equation}
where $\mathbf{W}$ is the weight matrix. 

The coefficient vector $\delta \mathbf{c}_{{k+1}_{l_2}}$ is known to include all possible coefficients from the dictionary of basis functions, resulting in a non-sparse representation. To address this, the research aims to achieve a minimal polynomial expansion that ensures sparsity in the generating function. Thus, a weighted $l_1$-norm optimization problem is proposed to minimize the number of coefficients selected from the extensive basis function dictionary. Unlike the equality constraint of Eq. \eqref{eq:FiniteDiff}, this optimization problem is formulated as a bounded two-norm error, incorporating $\epsilon$ as a soft inequality constraint. This approach allows sparse coefficients $\delta \mathbf{c}_{{k+1}_s}$ to balance sparsity and approximation error, offering increased flexibility. The complete optimization problem for selecting the minimal polynomial expansion of the generating function is outlined in Algorithm \ref{SparseSolution_algo}, with further details provided in Ref. \cite{jain2023stochastic,jain2023UP,jain2024stochastic}.

\begin{algorithm}[tb!] 
\caption{Iterative weighted $l_1$-norm optimization: $\delta\mathbf{c}_{{k+1}_s} $= WeightedOpt($\mathbf{A}, \mathbf{b},\mathbf{W},  \Delta_{s},\alpha, \epsilon , \eta  )$  \label{SparseSolution_algo}}
  \begin{algorithmic}[1] 
    \INPUT $\mathbf{A}, \mathbf{b},\mathbf{W},  \Delta_{s}, \alpha, \epsilon , \eta$
    \OUTPUT $\delta \mathbf{c}_{{k+1}_s}$
    \STATE \textbf{Initialization}   $\mathbf{K} \propto \mathcal{O}(basis) , \delta=1   $ 
    \STATE  compute $\delta\mathbf{c}_{k+1}^- =  \displaystyle \min _{\delta\mathbf{c}_{k+1}}\left\|\mathbf{K} \delta \mathbf{c}_{k+1}\right\|_{1} $  
    
     $\quad \qquad  \text{subject to:} \quad \left\|\mathbf{W}(\mathbf{A}  \delta \mathbf{c}_{k+1} - \mathbf{b} ) \right\|_{2}  \leq \epsilon$ 
    \WHILE{$\delta \ge \Delta_{s}$}
      \STATE Update $\mathbf{K}= \frac{1}{(\delta\mathbf{c}_{k+1}^- + \eta)}$,  
      $  \text{ 
          to find      }   \delta\mathbf{c}_{k+1}^+  =  \displaystyle \min _{\delta\mathbf{c}_{k+1}}\left\|\mathbf{K} \delta\mathbf{c}_{k+1}\right\|_{1} $  
      
      $\quad \quad \quad \quad \qquad \quad \qquad\quad \qquad \text{subject to:} \quad \left\|\mathbf{W}( \mathbf{A} \delta\mathbf{c}_{k+1} - \mathbf{b} ) \right\|_{2}  \leq \epsilon$
      \STATE  Compute $\delta = \| \delta\mathbf{c}_{k+1}^+ - \delta\mathbf{c}_{k+1}^-\|_2 $
      \STATE $\delta\mathbf{c}_{k+1}^- = \delta\mathbf{c}_{k+1}^+ $ 
    \ENDWHILE
    \STATE $\delta \mathbf{c}_{{k+1}_s} = \delta \mathbf{c}_{k+1}^-$
  \end{algorithmic}  
\end{algorithm} \setlength{\textfloatsep}{0pt}

This minimal representation $\delta \mathbf{c}_{{k+1}}$ is then employed to compute the generating function at time $t_{k+1}$, and the procedure can be repeated till the desired final time $t_f$. Now to obtain the solution of the trajectory tracking, given any initial perturbed state $\delta \x_0$ and the final boundary condition $\delta\boldsymbol{\psi}(\x(t_f))=\delta \x_f$, the initial co-state can be computed using coefficients at the final time ($\mathbf{c}_{t_f}$) using \eqref{eq:F2_x0} as:

\begin{equation}\label{eq:F2PHI_x0_tf}
\delta  \mathbf{x}_0=\frac{\partial F_2\left(\delta 
 \mathbf{x}, \delta  \boldsymbol{\lambda}_0, t_f\right)}{\partial (\delta \boldsymbol{\lambda}_0) } = \frac{\partial \Phi^T(\delta \mathbf{x}_f, \delta \boldsymbol{\lambda}_0) \mathbf{c}_{t_f} }{\partial (\delta \boldsymbol{\lambda}_0) }
\end{equation}
Once we obtain the value of the initial co-state corresponding to an initial and final perturbed state, the trajectory tracking problem is solved. The problem now becomes an initial value problem where we have the initial conditions on the perturbed state and co-state, and thus can be solved using the Non-linear MPC technique. A  similar procedure can also be applied to determine trajectory tracking from any perturbed state at the time instance $t_1$ to $t_f$. Now, let us look at some numerical examples to demonstrate the efficacy of the developed approach.

\section{Numerical Validation: PCR3BP}

\begin{table*}[t]
\caption{PCR3BP: Trajectory Tracking Problem}
\hspace{-0.15in}
\begin{tabular}{|c|c|}
\hline
Nominal OCP  &  Trajectory Tracking OCP \\
\hline \hline
\begin{minipage}{0.55\textwidth}
\begin{equation*} 
\begin{aligned}
\min_{\mathbf{u}(t)} \quad & J = \frac{1}{2} \int_{0}^{t_f} \mathbf{u}^T \mathbf{R} \mathbf{u} \, dt \\
& \text{subject to: } \dot{\mathbf{x}} = \mathbf{f}(\mathbf{x}, \mathbf{u}) \\
\mathbf{x}(t_0) = & [0.810796, -0.158270, -0.129473, 0.319169] \\
\mathbf{x}(t_{f}) = & [1.175974, -0.134272, -0.153277, -0.295254]
\end{aligned}
\end{equation*}
\end{minipage}
&
\begin{minipage}{0.45\textwidth}
\begin{equation*} 
\begin{aligned}
\min_{\mathbf{u}(t)} \quad & J = \frac{1}{2} \int_{0}^{t_f} \left(\delta \mathbf{x}^T \mathbf{Q} \delta \mathbf{x} + \delta \mathbf{u}^T \mathbf{R} \delta \mathbf{u}\right) dt \\
& \text{subject to: } \delta \dot{\mathbf{x}} = \mathbf{f}(\mathbf{x}, \mathbf{u}) - \dot{\mathbf{x}}_N \\
& \delta \mathbf{x}(t_0), \delta \mathbf{x}(t_{f}) \in \Omega
\end{aligned}
\end{equation*}
\end{minipage} \\ 
\hline \hline
Parameters & Test Case \\
\hline \hline
\begin{tabular}{c}
$\mathbf{R} = 1$, $\mathbf{Q} = 1$ \\
$t_f = 5 \text{ Days}$, $\mu = 0.0122$ \\
$\epsilon_s = 10^{-7}$, $\eta = 10^{-7}$ \\
$\delta_{rs} = 10^{-7}$, $\Delta_s = 10^{-8}$ \\
$\Omega(\delta\mathbf{x}_0) \in \pm [100 \text{ km}, 100 \text{ km}, 2 \text{ m/s}, 2 \text{ m/s}]$ \\
$\Omega(\delta\boldsymbol\lambda_0) \in \pm [0.02, 0.03, 0.005, 0.005]$ \\
$\text{Basis}(\mathcal{O}(4)) \xrightarrow{} m = 495$ \\
\\
\begin{tabular}{|c|}
\hline
Training and Testing Data \\
\hline \hline
$N_{train} = 1649 \in \Omega(\delta\mathbf{x}_0, \delta \boldsymbol\lambda_0)$ \\
$N_{test} = 100 \in \Omega(\delta \mathbf{x}_0, \delta \boldsymbol\lambda_0)$ \\
\hline
\end{tabular}
\end{tabular}
&
\begin{tabular}{c | c}
$\delta \mathbf{x}_0$ & $\delta \mathbf{x}_f$ \\[1.5mm]
\hline
\begin{tabular}{|c | c |}
\hline
Case & $[\delta x_0, \delta y_0]$ (km) \\
\hline
I & [100, 100] \\
II & [100, -100] \\
III & [-100, 100] \\
IV & [-100, -100] \\
\hline
\end{tabular} &
$\mathbf{0}$ \\[12mm]
\begin{tabular}{|c | c |}
\hline
Case & $[\delta \dot{x}_0, \delta \dot{y}_0]$ (m/s) \\
\hline
I-IV & $\mathcal{U} \in [-2, 2]$ \\
\hline
\end{tabular}
\end{tabular} \\
\hline \hline
\end{tabular}
\label{tab:CR3BP_OCP}
\end{table*}

\begin{figure*}
  \begin{subfigure}[b]{\columnwidth}
    \includegraphics[width=\linewidth]{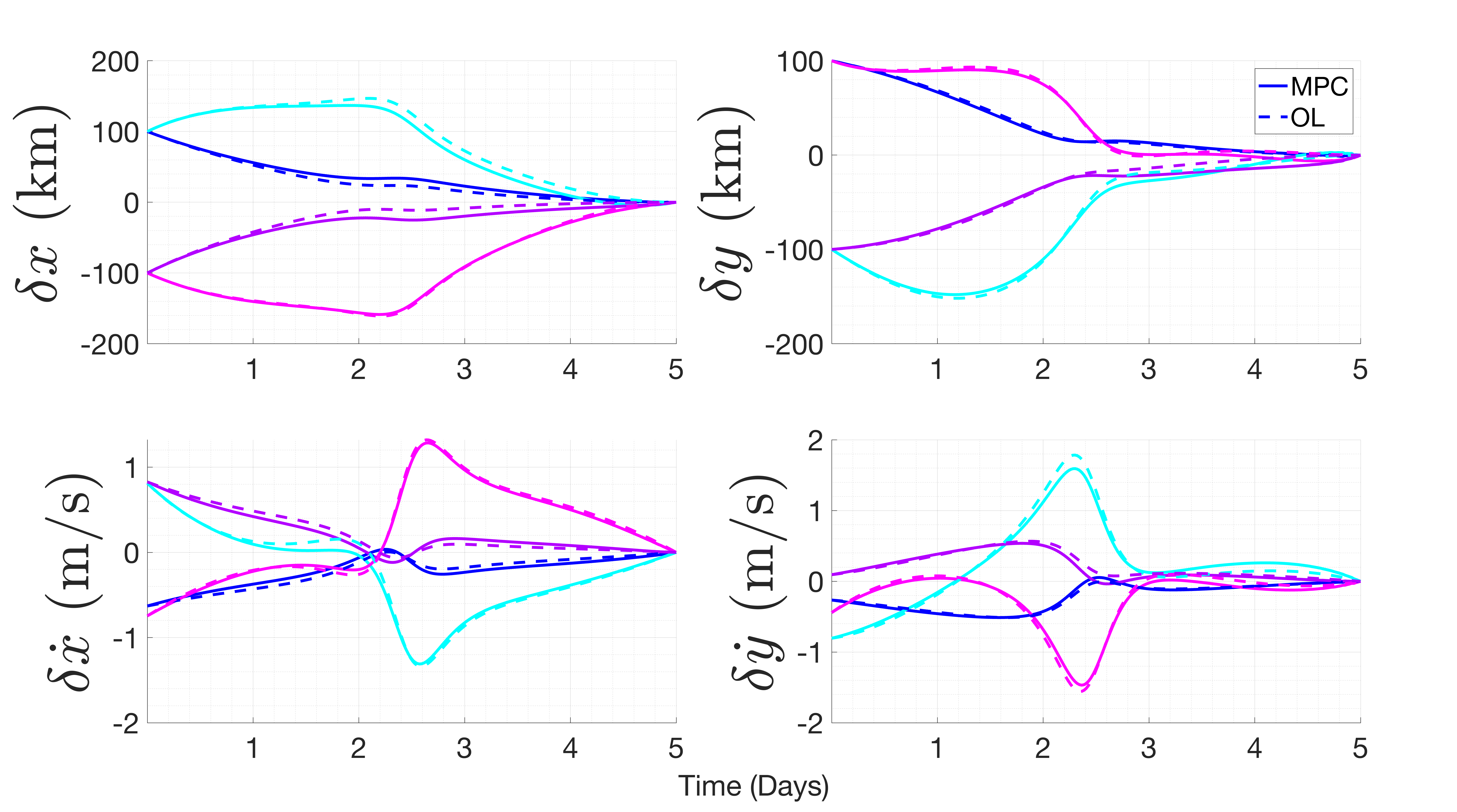}
    \caption{States Evolution}
    \label{fig:CR3BP_xy}
  \end{subfigure}
  \hfill 
  \begin{subfigure}[b]{\columnwidth}
    \includegraphics[width=\linewidth]{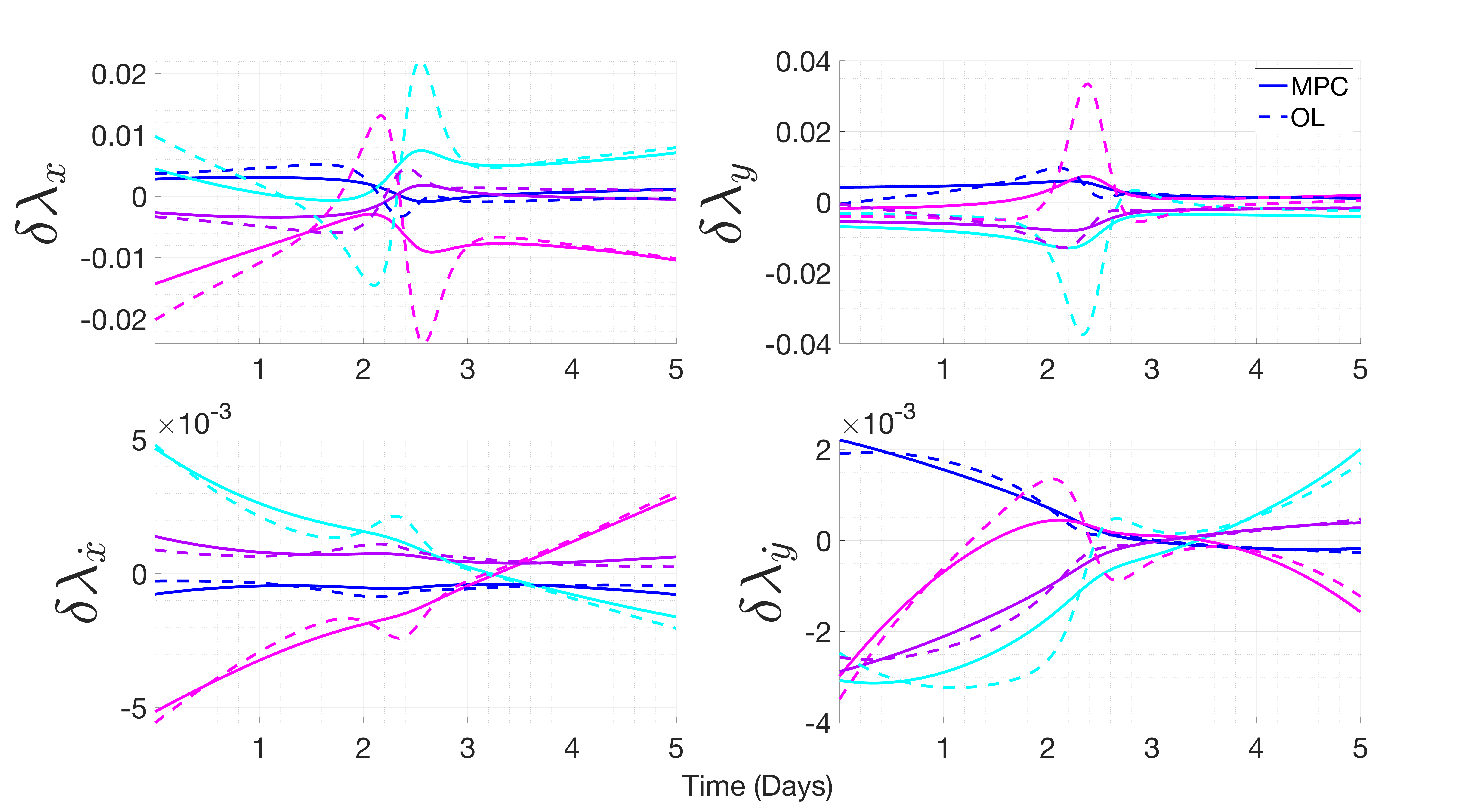}
    \caption{Costates Evolution}
    \label{fig:CR3BP_lambda}
  \end{subfigure}
  \caption{States and Costates evolution in the presence of initial perturbations}\label{fig:CR3BP_xlambda}
\end{figure*}

\begin{figure}
    \centering
    \includegraphics[width=\textwidth]{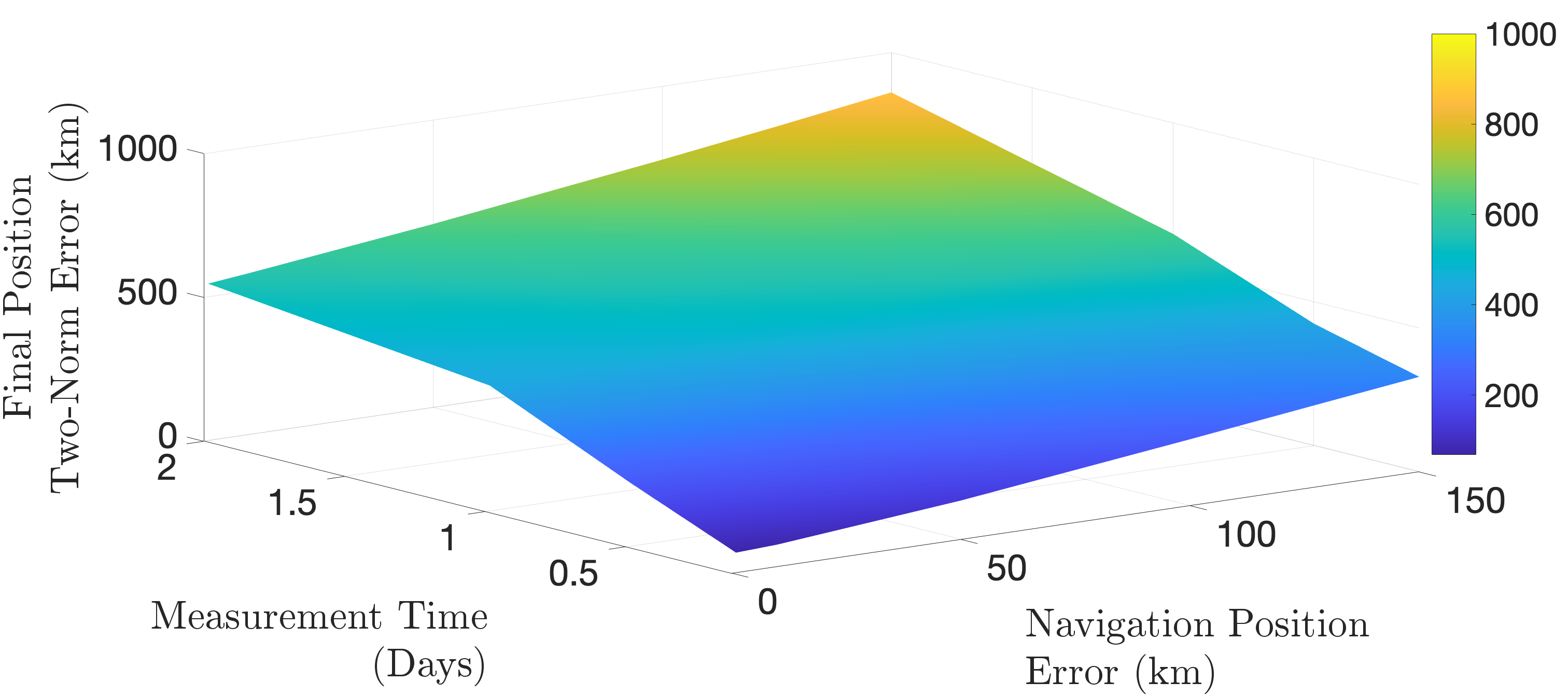}
    \caption{Variation of navigational position errors and the measurement time to reach the final position}
    \label{fig:Surf_NavStateError}
\end{figure}

This section applies the above methodology to track the transfer trajectory of a satellite from $L_1$ to $L_2$ Lyapunov orbit for a planar circular restricted three-body problem (PCR3BP) shown in \reffig{fig:CR3BP_nominal}. To calculate the nominal trajectory $\mathbf{x}_N$ from the $L_1$ to $L_2$ Lyapunov orbit, the cost function, along with the initial and final states in normalized coordinates, are detailed in Table \ref{tab:CR3BP_OCP}. The transfer time is specified as $t_f = 5$ days. Moreover, to track the nominal trajectory, the parameters, along with training and testing data, are presented in Table \ref{tab:CR3BP_OCP}. In this 8-dimensional system of states and costates, the training data comprises CUT samples totaling 1649, while for testing purposes, 100 random Monte Carlo samples are selected within the domain. Four distinct test cases are examined to demonstrate the efficacy of the developed method in accurately following the nominal trajectory. The initial state positions for these test cases are located at a corner of a hypercube with a side length of $200$ km, while the initial velocities are randomly selected from the range of $-2$ to $2$ m/s. 


The numerical results obtained from solving the Hamilton-Jacobi equation are shown in Table \ref{tab:CR3BP_OCP}, where the training and testing errors can be compared as a function of time. Notably, the training error consistently remains lower than the testing error across all time intervals. Additionally, the heatmap of the sparse coefficients is shown, revealing that basis functions up to third-order are actively engaged. In contrast, only a few fourth-order coefficients are involved towards the end.

To further validate the solution methodology, four initial perturbed states are selected to track the nominal trajectory using the developed approach. In \reffig{fig:CR3BP_xy}, the perturbed state variation for these four test cases is depicted along with the open-loop (OL) trajectories, demonstrating the successful tracking of the nominal trajectory and attainment of the final state within the specified time frame. Additionally, the time evolution of costates is presented in \reffig{fig:CR3BP_lambda}, showcasing that the perturbed trajectories effectively follow the nominal trajectory.

Furthermore, the fluctuations in the perturbation errors, coupled with navigational errors over time, are depicted in \reffig{fig:Surf_NavStateError}. In this figure, the velocity navigation error is assumed to vary randomly within $\pm 2 m/s$. Notably, when navigational uncertainties in position are minimal and measurements are taken frequently, the resultant final position two-norm errors are also minimal. Conversely, in scenarios where navigational position errors are substantial, compounded by infrequent measurements, the final position two-norm errors escalate significantly, impeding the achievement of the desired final position.


\FloatBarrier

\section{Summary}
This paper introduces a novel approach that leverages the Hamilton-Jacobi solution to improve the computation of optimal control in real time for tracking of paths in non-linear systems. By utilizing generating functions derived from the HJ equation, this method offers a robust solution that accommodates navigational uncertainties and model inaccuracies. The study showcases the potential of the Hamilton-Jacobi theory to mitigate the substantial computational costs associated with recomputing the optimal control solution at each stage of the MPC formulation. This is achieved by leveraging the generating function, which encapsulates all solutions within a single functional. The efficacy of the proposed methodology is demonstrated through MPC tracking on a PCR3BP transfer trajectory from $L_1$ to $L_2$.





\bibliographystyle{AAS_publication}
\bibliography{Biblio.bib,allrefs.bib}

\end{document}